\documentstyle[11pt,a4,twoside,amsmath,amssymb,amsfonts]{article}
\input xy
\xyoption{all}
\newtheorem{theorem}{Theorem}[section]
\newtheorem{lemma}[theorem]{Lemma}
\newtheorem{corollary}[theorem]{Corollary}
\newtheorem{proposition}[theorem]{Proposition}

\newcommand{\CH}{{\rm CH}}
\newcommand{\End}{{\rm End}}
\newcommand{\Ext}{{\rm Ext}}

\newcommand{\Hom}{{\rm Hom}}

\newcommand{\NS}{{\rm NS}}
\newcommand{\Pic}{{\rm Pic}}

\newcommand{\SL}{{\rm SL}}
\newcommand{\Quot}{{\rm Quot}}

\newcommand{\Tor}{{\rm Tor}}

\newcommand{\im}{{\rm im}}
\newcommand{\length}{{\rm length}}
\newcommand{\rk}{{\rm rk}}

\newcommand{\Ccal}{{\cal C}}

\newcommand{\Jcal}{{\cal J}}

\newcommand{\Ocal}{{\cal O}}

\newcommand{\cdop}{{\mathbb C}}

\newcommand{\pdop}{{\mathbb P}}
\newcommand{\zdop}{{\mathbb Z}}

\newcommand{\dual}{^\lor}
\newcommand{\ddual}{^{\lor \lor}}

\newcommand{\rar}{\rightarrow}

\newcommand{\rarpa}[1]{\stackrel{#1}{\rightarrow}}

\newcommand{\larpa}[1]{\stackrel{#1}{\leftarrow}}

\newcommand{\proof}{{\bf Proof:  }}

\newcommand{\bino}[2]{\left(
\begin{array}{c} #1 \\ #2 \end{array} \right)} 
\newcommand{\case}[1]{{\bf Case #1:}}
\newcommand{\step}[1]{{\bf Step #1:}}

\newcommand{\qed}{{ \hfill $\Box$}}

\newcommand{\tip}{{\vspace{0.5em} }}

\author{Georg Hein}
\title{Restriction of stable rank two vector bundles in arbitrary
characteristic}
\begin{document}
\maketitle
\begin{abstract}
Let $X$ be a  smooth variety defined over an algebraically closed field of
arbitrary characteristic and $\Ocal_X(H)$
be a very ample line bundle on $X$.
We show that for a semistable $X$-bundle $E$ of rank two,
there exists an integer $m$ depending only on $\Delta(E).H^{\dim(X)-2}$
and $H^{\dim(X)}$ such that the restriction of $E$ to a general divisor in
$|mH|$ is again semistable.
As corollaries we obtain boundedness results,
and weak versions of Bogomolov's theorem and Kodaira's vanishing theorem
for surfaces in arbitrary characteristic.
\end{abstract}

\section*{Introduction}

Let $(X,\Ocal_X(1)=\Ocal_X(H))$ be a smooth polarized variety
defined over an algebraic closed field
of arbitrary characteristic.
We assume $\Ocal_X(1)$ to be very ample.
Additionally, let $E$ be a $\mu$-semistable vector bundle of rank two on $X$.
We want to show that there exists an integer $m$ only depending on the
characteristic numbers $H^{\dim(X)}$ and
$(c_1(E)^2-4c_2(E)).H^{\dim(X)-2}$
such that the restriction of $E$ to a general element of $|mH|$ is
semistable.
Such effective bounds have been known only for the case that
the characteristic is zero.
In this case the restriction theorem of Flenner (see \cite{Fle}) gives
effective bounds on $m$ for semistable bundles of arbitrary rank.
On the other hand there are results of Mehta and Ramanathan
which say that the restriction of $E$ to a divisor in $|mH|$
is semistable (or stable, for $E$ a stable vector bundle) if $m \gg 0$ (cf.
\cite{MR1}, and \cite{MR2}).
A detailled overview on restriction theorems is given in \S7
of the book \cite{HL} of Huybrechts and Lehn.

First we discuss the case of rank two bundles on a surface $X$.
Theorem \ref{res1} shows that for a semistable $X$-vector 
bundle $E$ of rank two there exists an integer $m$ such that the
restriction of $E$ to
a general curve in the linear system $|mH|$ is semistable.

Using this result we provide a boundedness result for semistable rank two
bundles (proposition \ref{bound}).

For surfaces defined over $\cdop$ a semistable bundle $E$ cannot have
positive discriminant $\Delta(E)$
(Bogomolov's theorem, cf. \cite{Bog}).
In positive characteristic this does not hold.
No more than the Kodaira vanishing holds for positive characteristic (see
\cite{Ray}).
It is remarkable
that semistable bundles which contradict Bogomolov's theorem
behave well with respect to restrictions.
Applying our restriction result,
we obtain a weak form of Bogomolov's theorem (corollary \ref{res3},
cf.~also \cite{Meg} for vector bundles of arbitrary rank),
and a weak form of Kodaira's vanishing theorem in arbitrary characteristic
(corollary \ref{res4}).
The reader familiar with the vector bundle techniques presented in
Lazarsfeld's lectures, \cite{Laz} will deduce Reider type theorems for
surfaces in arbitrary characteristic.

If the $X$-bundle $E$ is semistable but not stable,
then it is easy to see that the restriction of $E$ to a general curve in
$|H|$ is semistable but not stable.
Conversely, we may ask whether stable bundles do restrict to stable
objects.
Theorem \ref{res5} gives an affirmative answer to this question.
The proof follows an idea of Bogomolov (see \cite{Bog1} and \cite{HL}
theorem 7.3.5) using the weak Bogomolov inequality deduced before.

Finally, we present with theorem \ref{res6} the higher dimensional version
of theorem \ref{res1}.
It turns out that its proof is easier than the proof in the surface case.
The main reason for this simplification is the fact that two general
hyperplanes in a linear system intersect in a irreducible subscheme. 

All these results should generalize to vector bundles of arbitrary rank.
To prove the corresponding results it seems necessary
to consider the complete Harder-Narasimhan filtration.

\tip
{\bf Acknowledgement:}
The author would like to thank D.~Huybrechts for many helpful
remarks.

\section{Preliminaries}
Let $X,H$ be a polarized projective variety.
We will identify line bundles on $X$ and their corresponding
Cartier divisor classes.
Moreover,
to any class in the Chow group $\CH^{\dim(X)}(X)$ of codimension
$\dim(X)$ cycles is assigned via evaluation on the fundamental class $[X]$
of $X$ its characteristic number.
This allows us to interpret $c_i(E).H^{\dim(X)-i}$ as integers.

For a coherent $X$-sheaf $E$,
we write $E(n)$ instead of $E \otimes \Ocal_X(H)^{\otimes n}$.
The Hilbert polynomial
$\chi_E:n \mapsto \chi(E(n))$ can be written in the 
following form
$$\chi_E(n) = a_0(E)\bino{n +\dim X}{\dim X}+
a_1(E)\bino{n +\dim X -1}{\dim X -1} + \ldots \,.$$
If $H$ is sufficiently general in the linear system $|H|$
(i.e., $\Tor^{\Ocal_X}_1(E,\Ocal_H)=0$),
we have $a_i(E)=a_i(E|_H)$,
for all integers $i < \dim X$.

We define the $H$-slope $\mu_H(E)$ of $E$
to be the quotient $a_1(E)/a_0(E)$.
A coherent sheaf $E$ is called Mumford semistable (resp. stable) with
respect to $H$,
if $E$ is torsion free,
and for all proper subsheaves $F \subset E$ the inequality
$\mu_H(F) \leq \mu_H(E)$ (resp. $\mu_H(F) < \mu_H(E)$)
holds true.
This kind of stability is also named slope stability, weak stability,
or $\mu$-stability.
For brevity we simply write stability
because we only use this stability concept.
We will frequently use the following facts on stable and semistable
coherent sheaves:
\begin{enumerate}
\item If $E$ and $F$ are semistable with $\mu_H(E)>\mu_H(F)$,
then the group $\Hom(E,F)$ vanishes.
\item For a stable bundle $E$ on a variety defined over an algebraically
closed field
the endomorphism group
$\End(E)$ consists of the scalar multiples of the identity.
\item If a rank two vector bundle $E$ is not semistable,
then there exists a unique maximal subsheaf $E_1 \subset E$ of rank one
which is the maximal destabilizing subsheaf.
Or equivalently,
there exists a unique destabilizing quotient $E \to Q$.
The flag $0 \subset E_1 \subset E$ is the Harder-Narasimhan filtration of
$E$.
\end{enumerate}

See, for example, the article \cite{Sha} of Shatz.

An important invariant of a vector bundle is its discriminant.
Let $E$ be a vector bundle of rank $r$ with Chern roots $\{\alpha_i
\}_{i=0,\ldots ,r}$.
As the name discriminant suggests we define the discriminant $\Delta(E)$
of the vector bundle $E$ by
$\Delta(E)= \sum_{i<j}(\alpha_i-\alpha_j)^2$.
Obviously $\Delta(E)$ can be expressed in terms of the Chern classes of
$E$, namely $\Delta(E)=(r-1)c_1(E)^2-2rc_2(E)$.
(Unfortunately, there are different definitions of $\Delta(E)$ in
literature,
differing by a sign or a constant.)
In particular, we have $\Delta(E)=c_1(E)^2-4c_2(E)$,
for a rank two vector bundle.

A rank two vector bundle $E$ on a surface $X$
is named Bogomolov unstable if there exists an injection
$A \rarpa{\iota} E$ of coherent sheaves where $A$ is an
$X$-line bundle, the cokernel of $\iota$ is torsion free,
and the inequalities  $(2A-c_1(E))^2>0$, and $(2A-c_1(E)).H>0$ are
satisfied for a polarization $H$ of $X$.

For a rational number $q$,
let $\lceil q \rceil$ be the least integer not smaller than $q$,
$\lfloor q \rfloor$ the largest integer smaller or equal to $q$, and
$[q]_+$ the maximum of $q$ and $0$.

\section{Rank two bundles on surfaces}
\subsection{The semistable restriction theorem}
\begin{theorem}\label{res1}
Let $X$ be a smooth surface over an algebraically closed field with a
very ample line bundle $\Ocal_X(1)=\Ocal_X(H)$.
For an $X$-vector bundle $E$ of rank two which is semistable with respect
to $\Ocal_X(1)$ the following holds:

\begin{enumerate}
\item If $\Delta(E) \geq 0$, then the restriction of $E$
to a general curve of the linear system $|H|$ is semistable;
\item For $\Delta(E)<0$ and any integer $l$ with
$l \geq \log_2 \left(\sqrt{\frac{-\Delta(E)}{H^2}}+1 \right)$
the restriction of $E$ to a general curve in $|2^lH|$ is semistable.
\end{enumerate}
\end{theorem}

\proof We divide the proof in several steps.
First we outline its strategy:
\begin{itemize}
\item We define the objects which are needed.
In particular, we define the non negative integer $A(m)$
which measures the instability of the restriction of $E$
to a general curve of the linear system $|mH|$
(step 1-3);
\item We next (step 4-7) compute an upper bound for $A(1)$.
This bound depends only on the Chern number $\Delta(E)$ and $H^2$;
\item After that, we give an upper bound for $A(2m)$ in terms of $A(m)$
(step 8-12);
\item Finally, we combine both estimates to conclude the theorem (step
13).
\end{itemize}

\step{1}
Let $m$ be a positive integer.
For the linear system $|mH|$ we denote by $\Ccal_{|mH|}$ the universal
curve over $|mH|$.
We have the morphisms
$$\xymatrix{|mH| & \Ccal_{|mH|} \ar[r]^q \ar[l]_p & X} .$$
The space $|mH|$ is isomorphic to $\pdop^{h^0(m)-1}$ where $h^0(m)$
denotes the dimension of $H^0(X,\Ocal_X(m))$.
Since $|mH|$ is base point free,
$q$ is a $\pdop^{h^0(n)-2}$-bundle.
We denote by $g_m$ the genus of a smooth curve of $|mH|$.
A curve $C\subset X$ rationally equivalent to $mH$ corresponds to a
geometric point in $|mH|$ which we denote by $[C]$.

\step{2}
For all integers $a$ with $2a<c_1(E).(mH)$ we consider the Quot scheme
$$\Quot_{m,a}:=\Quot^{P_a}_{q^*E / \Ccal_{|mH|} / |mH|}$$
of $p$-flat quotients of $q^*E$ with Hilbert polynomial
$P_a(k)=(mH^2)k+a+1-g_m$  with respect to the very ample line bundle
$q^*\Ocal_X(1)$ see \cite{Gro}.

For $a<g-1-h^1(E)-h^2(E(-mH))$ the scheme $\Quot_{m,a}$ is the empty
scheme.
To see this, we remark that for any curve $[C] \in |mH|$ we have the
inequality
$h^1(E|_C) \leq h^1(E)+h^2(E(-mH))$.
Hence, any quotient of the restriction $E|_C$ has at least Euler
characteristic $-(h^1(E)+h^2(E(-mH)))$.
From that bound, using the Riemann-Roch theorem for curves, we obtain the
above bound for the degree of quotients of $E|_C$.

Thus, we are considering only a finite number of Quot schemes.

\step{3}
Since the schemes $\Quot_{m,a}$ are projective over $|mH|$, 
their images dominate $|mH|$ if and only if at least one $\Quot_{m,a}$ is
surjective
over $|mH|$.
If they do not cover $|mH|$, then the restriction of $q^*E$ to the general
fiber of $p$ is semistable.
In this case we define the number $A(m)$ to be zero.
Otherwise we define $A(m)$ by
$$A(m):= \max \left\{c_1(E).(mH)-2a \left|
\begin{array}{ll}
2a<c_1(E).(mH) &\mbox{ and} \\
\Quot_{m,a} \to |mH| &\mbox{ is surjective.} 
\end{array}
\right. \right\} \,.$$
By definition $A(m)$ measures how far the restriction of $E$ to the
general curve of $|mH|$ is from being semistable.
The projectivity of the Quot schemes implies that if $A(m)>0$,
then the restriction of $E$ to any curve $[C]\in |mH|$ has a quotient $Q$
with Hilbert polynomial
$$\chi(Q(k))=(mH^2)k+\frac{1}{2}(c_1(E).(mH)-A(m))+1-g_m \, .$$
We will apply this specialization property
(see also \cite{Sha}) in the sequel
to reducible curves $[C] \in |2mH|$ with $C=C' \cup C''$
where $C'$ and $C''$ are smooth curves in $|mH|$,
to bound $A(2m)$ in terms of $A(m)$.

\step{4}
From now on we assume that $A(m)$ is positive.
We set $b(m):=\frac{1}{2}(c_1(E).(mH)-A(m))$.
By definition of $A(m)$ the subset
$$Y := \bigcup_{a<b(m)} \im( \Quot_{m,a} \to |mH|) $$
is a proper closed subset of $|mH|$.
For all points of the open subset $U'_m= |mH| \setminus Y$ the restriction
of $E$ the corresponding hyperplane has a minimal destabilizing quotient
of degree $b(m)$.
The minimality of the destabilizing quotient implies its uniqueness.
Therefore the restriction
$\Quot_{m,b(m)} \times_{|mH|} U'_m$ of the Quot scheme $\Quot_{m,b(m)}$ to
$U'_m$ gives a bijection of geometric points of
$\Quot_{U'_m}:=\Quot_{m,b(m)} \times_{|mH|} U'_m$ and $U'_m$.
Thus, $p_{U'_m}:\Quot_{U'_m} \to U'_m$ is an isomorphism or completely
inseparable.
If $[ \xymatrix{E|_C \ar@{->>}[r]^\alpha & F}]$ is a geometric point of
$\Quot_{U'_m}$,
then we have $\Hom(\ker (\alpha), F)=0$.
Therefore the relative tangent bundle of $p_{U'_m}$ vanishes.
Eventually, we conclude that $p_{U'_m}$ is an isomorphism.

\step{5}
If $[C]$ is a smooth curve in $U'_m$,
then the minimal quotient of degree $b(m)$ has to be a quotient line
bundle.
Therefore,
by considering the open subset $U_m$ of $U'_m$ parametrizing smooth curves,
we obtain the following situation:

$$\xymatrix{ U_m & \Ccal_{U_m} \ar[r]^q \ar[l]_p & X}$$
and a destabilizing quotient line bundle $L$ of $q^*E$ which is $p$-flat.
Furthermore the degree of $L$ on all fibers of $p$ is $b(m)$.
The surjection $q^*E \to L$ defines the following diagram:

$$\xymatrix{\Ccal_{U_m} \ar[d]_p \ar[r]^\xi \ar[dr]^q & \pdop(E) \ar[d]\\
U_m & X}$$
For a curve $C \subset X$ which is parametrized by $U_m$ we call
$\xi(p^{-1}[C])$ its canonical $m$-lifting.

\step{6}
Now we take two smooth curves $H_1$ and $H_2$ in $X$ which meet
transversally and which are contained in $U_1 \subset |H|$.
The pencil spanned by these curves defines a rational map
$\xymatrix{\pdop^1_{\mbox{ }} \ar@{-->}[r] & U_1 \ar[r]^-{\sim} & {\mbox{
}}} \!\!\Quot_{1,b(1)} \times U_1$.
Since $\Quot_{1,b(1)}$ is projective,
we obtain a morphism $\pdop^1 \to \Quot_{1,b(1)}$.
This corresponds to a flat family of degree $b(1)$ quotients for all
restrictions of $E$ to curves of the pencil.
To be precise we have the following situation:
$\pdop^1 \larpa{p} \tilde X \rarpa{q} X$
where $\tilde X$ denotes the blow up of $X$ in the points of $H_1 \cap
H_2$, and a destabilizing $p$-flat quotient $q^*E \to Q$ which for all
$p$-fibers is of degree $b(1)$.
Over $\pdop^1 \cap U_1$ the quotient $Q$ is a line bundle (see step 5).
$Q$ is flat and its restriction to most fibers is torsion free.
Hence, $Q$ itself is torsion free of projective dimension at most one.
The kernel $K$ of $q^*E \to Q$ is a line bundle on
$\tilde X$ which is isomorphic to $q^*L+\sum_{i=1}^{H^2} a_iE_i$
where the $\{ E_i \}_{i=1,\ldots,H^2}$ are the exceptional fibers of the
blow up $q$.
Restricting the exact sequence $0 \to K \to q^*E \to Q \to 0$
to $E_i$ we see that the integer $a_i$ is at least zero.

\step{7}
It results that $Q$ is of the form $\Ocal(q^*c_1(E)-q^*L-\sum a_i E_i)
\otimes \Jcal_Z$ where $\Jcal_Z$ denotes the ideal sheaf of a closed
subscheme $Z$ of $\tilde X$ of finite length.
Semistability of $E$ implies:
\begin{eqnarray}
H.D& \leq& 0 \qquad
\end{eqnarray}
where $D=(2L-c_1(E)) $.
Chern class computation gives:
$$\begin{array}{rcl}
c_2(E) & = & c_1(K)c_1(Q) +\length(Z) \\
& = & L.c_1(E)-L^2+\sum a_i^2 +\length(Z) \\
\end{array}$$
It follows
\begin{eqnarray}
\sum a_i^2 &\leq & c_2(E)+L^2-L.c_1(E) = c_2(E)-\frac{c_1(E)^2}{4}
+\frac{D^2}{4}
\end{eqnarray}
The discrepancy to semistability is the number $A(1)$
\begin{eqnarray}
A(1) & = & 2c_1(K).\left( H -\sum E_i\right) - c_1(E).H
=  D.H + 2\sum a_i
\end{eqnarray}
Now we use the inequality 
$$\sum_{i=1}^{H^2} a_i \leq \sqrt{H^2\sum_{i=1}^{H^2} a_i^2}$$
and inequalities (2) and (3) to deduce:
$$\begin{array}{rcl}
A(1) & \leq & D.H +2 \sqrt{H^2\cdot(c_2-\frac{c_1(E)^2}{4})+H^2 \cdot
\frac{D^2}{4}} \\
 & \leq & D.H +\sqrt{ -H^2\Delta(E) + H^2 \cdot D^2}
\end{array}$$

By the Hodge index theorem  $H^2 \cdot D^2 \leq (D.H)^2$.
Thus, we eventually obtain:
$$A(1) \leq D.H+\sqrt{(D.H)^2-H^2\Delta(E) } \, .$$
The basic properties of the function $x \mapsto x+\sqrt{x^2-H^2\Delta(E)}$
together with (1) give a bound for $A(1)$:
If $\Delta(E) \geq 0$, then $A(1)=0$.
For $\Delta(E) < 0$ we have the upper bound $A(1) < \sqrt{-H^2\Delta(E)}$.

\step{8}
Take a reducible curve $C=C' \cup C''$ where $C',C'' \in |mH|$ are smooth
curves which intersect transversally.
The singular divisor of $C$ consisting of  $m^2 \cdot H^2$
nodes we denote by $D$.
Let $E|_C \to Q$ be a torsion free quotient of $E$ with Hilbert polynomial
$$\chi_Q(k) = \chi( Q \otimes \Ocal_X(k))=2m \cdot H^2\cdot k+b+1-g_{2m}$$
where $g_{2m}$ denotes the arithmetic genus of $C$.
Torsion free means: $Q$ does not contain a subsheaf of dimension zero.
By the Mayer-Vietoris exact sequence
$$0 \to {\Ocal_C}\to {\Ocal_{C'} \oplus \Ocal_{C''}} \to {\Ocal_D} \to 0$$
we see that $g_{2m}=2g_m+ m^2 \cdot H^2-1$.
Furthermore,
we obtain from this exact sequence the following diagram
with exact rows and surjective columns
$$\xymatrix{0 \ar[r] & E|_C \ar[d] \ar[r]
& E|_{C'} \oplus E|_{C''} \ar[d] \ar[r] & E|_D \ar[d] \ar[r] & 0 \\
\Tor_1^{\Ocal_C}(Q,\Ocal_D) \ar[r]  & Q  \ar[r]
& Q|_{C'} \oplus Q|_{C''}  \ar[r] & Q|_D  \ar[r] & 0}$$
Since $\Tor_1^{\Ocal_C}(Q,\Ocal_D)$ is concentrated in $D$,
and $Q$ was assumed to be torsion free,
the image of $\Tor_1^{\Ocal_C}(Q,\Ocal_D)$ in $Q$ is zero.
Therefore, the equality
\begin{eqnarray}\label{chiq}
\chi(Q)+\length(Q|_D)& =&\chi( Q|_{C'}) +\chi( Q|_{C''})
\end{eqnarray}
holds true.
There are three cases for the ranks of $Q|_{C'}$ and $Q|_{C''}$.
The pair $(\rk(Q|_{C'}),\rk(Q|_{C''}))$ has to be
$(1,1)$, $(2,0)$, or $(0,2)$.

\step{9}
We next show that if the ranks of $Q|_{C'}$ and $Q|_{C''}$ do not coincide,
then the quotient $Q$ is not destabilizing.
Assume that $Q|_{C'}$ has rank two and $Q|_{C''}$ is torsion.
It results that $Q|_{C'}$ is isomorphic to $E|_{C'}$, and $Q|_{C''}$ is
isomorphic to $E|_D$.
Hence, $Q$ is isomorphic to $E|_{C'}$.
Therefore we find
$$\chi(Q(k))=\chi(E(k))-\chi(E(k-m))=(2mH^2)k+\chi(E)-\chi(E(-m)) \,.$$
Analogously we compute the Euler characteristic of $E|_C$ to be
$$\chi(E|_C(k))=(4m \cdot H^2)k+\chi(E)-\chi(E(-2m)) \,.$$
In order to prove that $Q$ is not destabilizing we must show
that the inequality
$$\frac{\chi(E)-\chi(E(-m))}{2mH^2} >
\frac{\chi(E)-\chi(E(-2m))}{4mH^2}$$
holds.
This inequality is equivalent to 

$$2(\chi(E)-\chi(E(-m))) >
\chi(E)-\chi(E(-2m)) \,.$$
The last inequality holds because the function $k \mapsto \chi(E(k))$ is
strictly convex by the Riemann-Roch theorem for surfaces.

\step{10}
{\bf (General intersection lemma) \/}{\em
Let $C''$ be an irreducible curve in $X$ with a lifting $\tilde C''$ to
$\pdop(E)$.
For a general curve $[C']$ in $U_m$ its canonical $m$-lifting $\tilde
{C'}$ in $\pdop(E)$ intersects $\tilde C''$ in zero or $C'.C''$ points.}

\proof For brevity, we write  $U$ instead of $U_m$.
We consider the following situation:
$$\xymatrix{& \Ccal_U \times_X C'' \ar[r] \ar[d] \ar[ddr]^(0.3)\psi
& C'' \ar[d]^\iota \\
\Ccal_U \times_{\pdop(E)} \tilde C'' \ar[ur]^\varphi \ar[r] \ar[d]
& \Ccal_U \ar[r]|\hole \ar[d] \ar[dr]|\hole & X \\
\tilde C'' \ar[r]^{\tilde \iota} & \pdop(E) \ar[ur]|(0.65)\hole & U}$$
Since $C'' \to X$ and $\Ccal_U \to U$ are projective morphisms,
so is $\psi$.
The same way, we see that the composition morphism $\psi \circ \varphi$ is
projective.
For a geometric point $[C'] \in U$ the fiber of $\psi$ over $[C']$ is the
intersection $C' \cap C''$.
Thus, $\psi$ is of relative dimension zero.
Analogously we identify the fiber of $\psi \circ \varphi$ with the
intersection of the liftings.
By construction $\Ccal_U$ is an open subset in a $\pdop^n$-bundle over 
$X$.
We conclude the irreducibility of $\Ccal_U \times_X C''$.
Thus we see,
that if $\varphi$ is a dominant morphism,
then the canonical lifting $\tilde C'$ of a general curve $C'$  intersects
$\tilde C''$ in $C'.C''$ points.
If the morphism $\varphi$ is not dominant,
then the canonical lifting of a general curve $C'$ is disjoint from
$\tilde C''$.
\hfill{$\Box$ (of step 10)}

\step{11}
Let $[C'']$ be a point in $U_m$, and $\tilde C''$ be its canonical lifting
to $\pdop(E)$.
Let us assume that the lifting $\tilde C'$ of a general curve $[C'] \in
U_m$ intersects $\tilde C''$ in $C'.C''$ points.
If we consider the pencil spanned by $C''$ and $C'$,
then we obtain (see step 6) a family $Q$ over this pencil with all
$a_i$ equal to zero.
Indeed, if one $a_i$ is positive,
then the $m$-lifting of general curve contained in the pencil spanned by
$C'$ and $C''$ does not intersect $\tilde C''$ in the point $P_i$.
This would imply (see (1) and (3) of step 7) that $A(m)=0$.

By the above lemma we can assume that
$[C'],[C''] \in U_m$ are two smooth curves
whose canonical liftings are disjoint.
We consider now for $C=C' \cup C''$ a minimal quotient $Q$ of $E|_C$
having rank one on $C'$ and $C''$.
We call a point $P \in D=C' \cap C''$ a point of discord if the dimension
of $Q \otimes k(P)$ is two.
Let $M$ be the number of points of discord.
It is obvious that the maximal torsion subsheaf of $Q|_{C'}$
is concentrated in the points of discord and has length $M$.
The quotient of $Q_{C'}$ modulo its torsion is denoted by $Q'$,
and analogous we have the $C''$-line bundle $Q''$.
We obtain from (\ref{chiq}) that
\begin{eqnarray}\label{chiq2}
\chi(Q) & = & \chi(Q') + \chi(Q'') + M -m^2H^2
\end{eqnarray}

\step{12} {\em The inequality $A(2m) \leq [2A(m)-2m^2H^2]_+$
holds.}

\proof
We consider the unique destabilizing quotient $L'$ of $E|_{C'}$.
The kernel of $ E|_{C'} \to L$ we denote by $F'$.
We now consider the composition $\beta':F' \to E|_{C'} \to Q'$.

\case{1} The morphism $\beta'$ (or $\beta''$) is not trivial.

If $\beta'$ is not trivial, then it follows that $\deg(Q') > \deg(F')$ and
$\chi(Q'(k)) \geq \chi(F'(k)) = \chi(F''(k))$.
We obtain from (\ref{chiq2}) that
$$\begin{array}{rcl}
\chi(Q(k)) & \geq & \chi(F''(k))+ \chi(Q''(k))-m^2H^2 \\
& \geq &  \chi(F''(k))+ \chi(L''(k))-m^2H^2 \\
& = & \chi(E|_{C''}(k))-m^2H^2 \\
& = & \frac{1}{2}(\chi(E|_C(k))) \,.
\end{array}$$
Thus, in this case the quotient $Q$ is not destabilizing.

\case{2} If $\beta'$ is trivial, then we obtain
$F' \subset \ker(E_{C'} \to Q')$.
However $F'$ is the rank one subbundle of $E_{C'}$
of maximal degree. Hence we have $Q' \cong L'$.
Since the canonical liftings of $C'$ and $C''$ do not intersect,
we must have $M=m^2H^2$.
The equality (\ref{chiq2}) consequently yields:
$$\begin{array}{rcl}
\chi(Q(k)) & = & 2\chi(L'(k))\\
&=&\chi(E|_C'(k))-A(m) \\
&=& \frac{1}{2}\chi(E|_C(k))+(m^2H^2-A(m))
\end{array}$$
This gives the asserted inequality for $A(2m)$.

To complete the proof of step 12 we just remark
that a similar computation shows that $A(m)=0$ implies $A(2m)=0$.
\hfill{$\Box$ (of step 12)}

\step{13}
Using induction on $l$,
we obtain from the inequality of step 12 
$$A(2^l) \leq [ 2^lA(1)-2^l(2^l-1)H^2]_+ \,.$$
Combining this with the upper bound for $A(1)$ computed in step 7 the
theorem follows.
\qed

\tip
{\bf Remark 1:}
The theorem still holds true for a semistable coherent $X$-sheaf $E$ of
rank two.
Indeed, consider the embedding $E \rar E\ddual$ of $E$ into its double
dual.
In this way, we obtain the semistable vector bundle $E\ddual$ which
obviously
satisfies $\Delta(E\ddual) \geq \Delta(E)$.
Hence, the theorem applies.

\tip
{\bf Remark 2:}
We can extend the theorem to projective surfaces $X$ with isolated
singularities.
Reviewing the proof,
we see that it is enough to have smooth curves in the linear systems
$|mH|$.
By the same argument, we see that the theorem holds true
if we require
$\Ocal_X(1)$ to be base point free.

\tip
{\bf Remark 3:}
If we know the ideal
$\{ L.H \}_{L \in \Pic(X)} \subset \zdop$
of intersections with $H$,
then we can sharpen the inequality of step 7.
To illustrate this,
let us assume that the Picard group of $X$ is generated by $\Ocal_X(1)$.
Furthermore,
suppose that $E$ is a semistable $X$-vector bundle.
We have $\det(E)=nH$.
If $\Delta(E)<0$, then we can improve the bound for $A(1)$ of step 7 by
$$A(1) \leq \left\{ \begin{array}{ll}
-2H^2+\sqrt{4(H^2)^2-H^2\Delta(E)} & \mbox{ for } n \mbox{ even;} \\
-H^2+\sqrt{(H^2)^2-H^2\Delta(E)} & \mbox{ for } n \mbox{ odd.} \\
\end{array} \right.$$

\subsection{A boundedness result}
Let $X$, $\Ocal_X(1)$ be as before.
Furthermore, let $E$ be a semistable $X$-vector bundle of rank two.
We next give a bound $M$ depending only on the characteristic numbers
$c_2(E)$ and $c_1(E).H$ such $E(M)$ becomes globally generated.
We use Mumford's concept of $m$-regularity:

A coherent sheaf $E$ on a polarized variety $X$ with very ample line
bundle $\Ocal_X(1)$ is called $m$-regular, if $h^i(E(m-i))=0$, for all
$i>0$.

The following lemma (cf. \S 14 in \cite{Mum}) resumes properties of
$m$-regular sheaves.

\begin{lemma}\label{bound1}
Let $X$ be a projective variety with a very ample line bundle
$\Ocal_X(1)=\Ocal_X(H)$, and $E$ be a coherent $X$-sheaf.
If $E$ is $m$-regular,
then it is globally generated, and
$E(m+k)$-regular for all $k \geq 0$.

Let $D \in |H|$ be a divisor such that the sequence
$0 \to E(-1) \to E \to E|_D \to 0$ is exact.
If $E|_D$ is $m$-regular, then $E$ is $(m+h^1(E(m-1)))$-regular.
\qed
\end{lemma}

This lemma outlines our strategy.
We first show that for a suitable curve $C \in |H|$ and an integer $m_1$ the
restriction $E_C=E|_C$ is $(m_1+1)$-regular.
In order to obtain the boundedness result,
we then compute an upper bound for $h^1(E(m_1))$.

\begin{lemma}\label{bound2}
Let $E_C$ be a rank two vector bundle on a smooth curve $C$ of genus $g$
defined over an algebraically closed field.
We define the number $A$ to be zero if $E_C$ is semistable.
Otherwise we set
$$A = \max \{ \deg(E_C)-2\deg(Q) \, | \,
E_C \to Q \mbox{ is surjective, and }\rk(Q)=1 \} \, .$$
\begin{enumerate}
\item[(1)] If $L$ is a $C$-line bundle with
$\deg(L)>\frac{A-\deg(E_C)}{2}+2g-2$,\\
then $H^1(C,E_C \otimes L)=0$;
\item[(2)] For any $C$-line bundle $L$ the inequality\\
$h^0(E_C\otimes L) \leq 2 \left[ 1+\deg(L)-\frac{\deg(E_C)-A}{2}
\right]_+$ holds true.
\end{enumerate}
\end{lemma}
\proof
(1) If $h^1(E_C \otimes L)>0$, then there exists,
by Serre duality,
a non trivial homomorphism $\varphi:E \to \omega_C \otimes L^{-1}$.
(Here $\omega_C$ denotes the dualizing sheaf of $C$.)
Thus,
the image of $\varphi$ is a rank one quotient of degree at most
$2g-2-\deg(L)$.
By the very definition of the number $A$ we obtain
$\deg(L) \leq \frac{A-\deg(E_C)}{2}+2g-2$.

(2) Analogously, we see that for $\deg(L)<\frac{\deg(E_C)-A}{2}$ there are
no global sections of $E_C \otimes L$.
Thus, the assertion holds for all line bundles $L$ of degree less than
$\frac{\deg(E_C)-A}{2}$.
Let $P \in C$ be a geometric point of $C$.
Then from the exact sequence
$0 \to E_C \otimes L(-P) \to E_C \otimes L \to E_C|_P \to 0$
we obtain $h^0(E_C \otimes L) \leq 2+ h^0(E_C \otimes L(-P))$
which proves the second statement.
\qed

\tip
We now take a smooth curve $C$ of genus $g$ in the linear system $|H|$
such that for the restriction $E_C$ the number $A$ of the above lemma is
at most $\sqrt{ [-H^2 \cdot \Delta(E)]_+}$.
We have seen in step 7 of the proof of theorem \ref{res1}
that this is possible.
The adjunction formula gives $2g-2=H.(H+K_X)$.
Obviously, the degree of the $C$-line bundle $\Ocal_X(mH)|_C$ is $mH^2$.
Thus, setting
$$m_1:= \left\lfloor \frac{1}{H^2} \left( \frac{\sqrt{ [-H^2 \cdot
\Delta(E)]_+}-c_1(E).H}{2}+H.(H+K_X) \right) \right\rfloor +1$$
we obtain by lemma \ref{bound2} that $E_C$ is $(m_1+1)$-regular. 

The semistability of $E$ implies that $h^0(E(m_2-1))=0$, for
$m_2:=\left\lceil \frac{-H.c_1(E)}{2H^2} \right\rceil$.
Applying the inequality $h^0(E(m))<h^0(E(m-1))+h^0(E_C(m))$
obtained from the long exact cohomology sequence yields
$$h^0(E(m_1)) \leq m_3 := 2\sum_{m=m_2}^{m_1}
\left[ 1+ mH^2-\frac{c_1(E).H-\sqrt{ [-H^2 \cdot \Delta(E)]_+}}{2}
\right]_+ \, .$$
Since $h^2(E(m_1))=0$ we deduce that $h^1(E(m_1)) \leq m_3-\chi(E(m_1))$.
Setting $m_4:=m_1+m_3-\chi(E(m_1))$, we obtain by lemma \ref{bound1}:

\begin{proposition}\label{bound}
Let $X$ be a smooth projective surface over an algebraically closed field,
and $\Ocal_X(1)$ a very ample line bundle on $X$.
Furthermore, let $E$ be a rank two $X$-bundle which is semistable with
respect to $\Ocal_X(1)$.
Then for $m \geq m_4$ we have that $E(m)$ is globally generated.
The number $m_4$ defined above depends only on the characteristic numbers
of $E$.
\qed
\end{proposition}

It follows that any semistable sheaf $E$ of rank two with given $c_1(E).H$,
$c_1(E).K_X$, and $c_2(E)$ is a quotient of $\Ocal_X(-m_4)^{\oplus
\chi_E(m_4)}$.
Considering the Quot scheme $\Quot^{\chi_E}_{\Ocal_X(-m_4)^{\oplus
\chi_E(m_4)}/X}$
together with its natural $\SL_{\chi_E(m_4)}$-action
we obtain (see \cite{GIT}) the coarse moduli space of semistable coherent
sheaves on $X$ with Hilbert polynomial $\chi_E$.
This proves the next corollary.

\begin{corollary}
There exists a projective coarse moduli space for semistable
coherent sheaves of rank two with fixed
characteristic numbers on a smooth projective surface.
\end{corollary}

\subsection{Further applications}

\begin{proposition}\label{res2}
Let $X$ be a smooth projective surface
over an algebraically closed field
with a very ample line bundle $\Ocal_X(1)=\Ocal_X(H)$.
If $E$ is a $X$-bundle,
which is stable with respect to $\Ocal_X(1)$ and of
rank 2, then the inequality
$$\Delta(E) \leq \left\{ \begin{array}{ll}
1-4\chi(\Ocal_X) & \mbox{if } K_X.H<0 \,; \\
2-4\chi(\Ocal_X) & \mbox{if } K_X.H=0 \,; \\
\left[ 6-4\chi(\Ocal_X) +4\cdot \left\lceil \frac{K_X.H}{H^2} \right\rceil
K_X.H \right]_+ & \mbox{if } K_X.H>0 \,. \\
\end{array} \right.$$
holds.
\end{proposition}
\proof
We compute, using the Riemann-Roch theorem for surfaces
that
$$\chi(E \otimes E\dual) = \Delta(E)+4\chi(\Ocal_X)\,. $$
The stability of $E$ implies that
$H^0(E \otimes E\dual) = \Hom(E,E)$ is of dimension one.
Now we want to bound $h^2:=h^2(E \otimes E\dual)$.
By Serre duality,
$h^2$ equals the dimension of $\Hom(E,E(K_X))$
where $K_X$ denotes the canonical class on $X$.
Thus, we obtain for $K_X.H \leq 0$ 
$$h^2(E \otimes E\dual) \leq \left\{
\begin{array}{ll}
0 & \mbox{if } K_X.H < 0 \, ;\\
1 & \mbox{if } K_X.H = 0 \, .\\
\end{array} \right.$$
If $K_X.H >0$,
we set $m=\lceil \frac{K_X.H}{H^2}\rceil$
and consider a smooth curve $C$ in the linear system $|mH|$.
If $\Delta(E) \geq 0$, then by theorem \ref{res1},
we may assume that the restriction $E|_C$ is semistable.
Thus, $\Hom(E|_C,E|_C)$ is at most of dimension 4.
By induction we see
that for a $C$-line bundle $L$ of degree $d$ we can bound the
dimension of $\Hom(E|_C,E|_C \otimes L)$, by $4+4\cdot d$.

By definition of $m$,
we have $(K_X-mH).H \leq 0$.
Thus, we can bound the dimension of $\Hom(E,E(K_X-C))$ by one.
From the exact sequence
$$0 \to \Hom(E,E(K_X-C)) \to \Hom(E,E(K_X)) \to \Hom(E,E(K_X)|_C) $$
we obtain the estimate
$$h^2(E \otimes E\dual) \leq 5+4 \cdot
\left\lceil \frac{K_X.H}{H^2}\right\rceil K_X.H \, .$$
Applying the obvious inequality $\chi(E \otimes E\dual) \leq h^0(E \otimes
E\dual) + h^2(E \otimes E\dual)$ we obtain the estimation of the
proposition.
\qed

\begin{corollary}\label{res3}
{\bf (Weak Bogomolov inequality)}
Let $X,H$ be a very ample polarized smooth surface over an algebraically
closed field.
Let $E$ be a rank 2 vector bundle on $X$ satisfying

$$\Delta(E) > \left\{ \begin{array}{ll}
[1-4\chi(\Ocal_X)]_+ & \mbox{if } K_X.H<0 \,; \\
\left[ 2-4\chi(\Ocal_X)\right]_+ & \mbox{if } K_X.H=0 \,; \\
\left[ 6-4\chi(\Ocal_X) +4\cdot \left\lceil \frac{K_X.H}{H^2} \right\rceil
K_X.H \right]_+
 & \mbox{if } K_X.H>0 \,. \\
\end{array} \right.$$
Then $E$ is Bogomolov unstable.
\end{corollary}
\setcounter{equation}{0}
\proof By proposition \ref{res2} $E$ cannot be stable with respect to the
given polarization $H$.
Thus, we have a short exact sequence
$$0 \to A \to E \to \Jcal_Z(c_1(E)-A) \to 0 \, ,$$
where $Z\subset X$ is a closed subscheme of codimension 2.
Since $A$ is destabilizing we have 
$(c_1(E)-2A).H  \leq  0$.
Using the exact sequence to compute $c_2(E)$ yields
$$(c_1(E)-2A)^2 = \Delta(E)+4\cdot \length(Z) > 0 \,.$$
Thus,
the Hodge index theorem implies $(c_1(E)-2A).H<0$.
Consequently, $E$ is not semistable.
Now the half c\^one in the N\'eron-Severi group $\NS(X)$
defined by positive self intersection and negative
intersection with an ample class $H$ does not depend on $H$.
\qed

\tip
{\bf Remark:} G. Megyesi proves that for a vector bundle $E$ of arbitrary
rank on a smooth surface defined over a field of characteristic $p>0$ with
$\Delta(E)>0$ the pullback $(F^n)^*E$ of $E$ by a large power $n$ of the
absolute Frobenius $F$ is Bogomolov unstable (see \cite{Meg}).
Corollary \ref{res3} gives an effective bound for $n$,
for a vector bundle $E$ of rank two.

\begin{corollary}\label{res4}
{\bf (Weak Kodaira vanishing)}
Let $X$ be a smooth projective surface
defined over an algebraically closed field
with very ample line bundle $\Ocal_X(H)$.
Let $L$ be a nef $X$-line bundle such that
$$L^2> \left\{ \begin{array}{ll}
\left[ 1-4\chi(\Ocal_X) \right]_+ & \mbox{if } K_X.H<0 \,; \\
\left[ 2-4\chi(\Ocal_X)\right]_+ & \mbox{if } K_X.H=0 \,; \\
\left[ 6-4\chi(\Ocal_X) +4\cdot \left\lceil \frac{K_X.H}{H^2} \right\rceil
K_X.H \right]_+
 & \mbox{if } K_X.H>0 \,. \\
\end{array} \right.$$
Then the first cohomology group $H^1(X,L^{-1})$ vanishes.
\end{corollary}
\proof
\setcounter{equation}{0}
We take an extension $E$ of $\Ocal_X$ by $L^{-1}$.
Since $H^1(X,L^{-1}) = \Ext^1(\Ocal_X, L^{-1})$,
we have to show that the short exact sequence
$$0 \to L^{-1} \to E \to \Ocal_X \to 0$$
splits.
We compute $c_1(E)=-L$, $c_2(E)=0$, and $\Delta(E)=L^2$.
Consequently, by corollary \ref{res3},
$E$ has a destabilizing subsheaf $A$ of rank one with
$(2A+L)^2>0$, and
$(2A+L).H > 0$, for all ample classes $H$.
Since nef bundles are limits of ample classes,
we obtain
\begin{eqnarray}\label{r41}
(2A+L).L & \geq &0 \,.
\end{eqnarray}
By the same reason, the Hodge index theorem applies
\begin{eqnarray}\label{r42}
A^2 & \leq & \frac{(A.L)^2}{L^2} \,.
\end{eqnarray}
The subsheaf $A$ of $E$ cannot be contained in $L^{-1}$
because $A$ destabilizes $E$ whereas $L^{-1}$ does not.
Thus, $A$ is contained in $\Ocal_X$.
We conclude
\begin{eqnarray}\label{r43}
A.L & \leq & 0 \,.
\end{eqnarray}
Computing the second Chern class of $E$ in terms of $A$ and
$E/A$ we obtain
\begin{eqnarray}\label{r44}
A.L+A^2 \geq & 0
\end{eqnarray}
Combining (\ref{r42}) and (\ref{r44}) yields $A.L \leq -L^2$
or $A.L \geq 0$.
In view of (\ref{r41}) and (\ref{r43}) we deduce that $A.L=0$.
This equality, and the inequalities (\ref{r42}) and (\ref{r44})
imply $A^2=0$.

Now we claim that $A.H=0$.
Suppose this were not the case.
Then $(L+a\cdot A).H$ would be zero,
for a rational number $a$.
Applying once again the Hodge index theorem yields
$(L+a\cdot A)^2=L^2 \leq 0$
which contradicts our assumptions on $L$.

Since $A$ is contained in $\Ocal_X$ we conclude that $A=\Ocal_X$.
Thus the injection $A \to E$ splits the exact sequence.
\qed

\subsection{The stable restriction theorem}
\begin{theorem}\label{res5}
Let $X$ be a smooth projective surface over an algebraically
closed field with very ample line bundle $\Ocal_X(H)$.
Let $E$ be an $X$-vector bundle of rank 2 which is stable
with respect to the polarization $H$.
Assume furthermore, that the positive integer $m$ satisfies
\begin{itemize}
\item $m>\frac{a}{2H^2}-\frac{\Delta(E)}{2a}$
where $a>0$ is an integer not larger than the
positive generator of the ideal $\{ A.H \, | \, A \in \Pic(X) \}$;
\item $m^2H^2+\Delta(E)>\left\{ \begin{array}{ll}
\left[ 1-4\chi(\Ocal_X) \right]_+
& \mbox{if } K_X.H<0 \,; \\
\left[ 2-4\chi(\Ocal_X) \right]_+
& \mbox{if } K_X.H=0 \,; \\
\left[ 6-4\chi(\Ocal_X) +4\cdot \left\lceil \frac{K_X.H}{H^2} \right\rceil
K_X.H \right]_+
 & \mbox{if } K_X.H>0 \,. \\
\end{array} \right.$
\end{itemize}
Then the restriction of $E$ to any smooth curve of the
linear system $|mH|$ is stable.
\end{theorem}
\proof
\setcounter{equation}{0}
Let $C$ be a smooth curve rationally equivalent to $mH$.
Suppose that $E|_C$ has a quotient line bundle $L$ of degree $d$
where
\begin{eqnarray}\label{r31}
2d & \leq & c_1(E).(mH) \, .
\end{eqnarray}
We denote the kernel of the surjection $E \to L$ by $E'$.
Then:
\begin{eqnarray}\label{r32}
c_1(E')=c_1(E)-mH & &  c_2(E')=c_2(E)+d-c_1(E).(mH)
\end{eqnarray}
Hence, we obtain from (\ref{r31}) 
$$
\Delta(E') = \Delta(E)+m^2H^2+2(c_1(E).(mH)-2d) 
\geq  \Delta(E)+m^2H^2 \, .$$
From proposition \ref{res2} and the assumptions on $m$,
we conclude that $E'$ cannot be stable with respect to $H$.
Thus there is a line bundle $A$ destabilizing $E'$.
We may assume that $E'/A$ is torsion free.
We consider the following short exact sequence.
$$0 \to A \to E' \to \Jcal_Z(c_1(E)-mH-A) \to 0 $$
where $Z\subset X$ denotes a closed subscheme of codimension 2.
Since $A$ is destabilizing we obtain
$H.(2A+mH-c_1(E)) \geq 0$.
We want to show that $A$ destabilizes $E$.
It is convenient to introduce the divisor $B:=2A-c_1(E)$.
Therefore we have to show that $B.H \geq 0$.
Using this notation the last inequality reads
\begin{eqnarray}\label{r33}
B.H & \geq & -mH^2 \,.
\end{eqnarray}

Computing $c_2(E)$ via the above exact sequence yields
$$c_2(E')=A.(c_1(E)-mH-A)+\length(Z) \,.$$
Using $\length(Z) \geq 0$, (\ref{r31}), and (\ref{r32}) we obtain
$$A.(c_1(E)-mH-A) \leq c_2(E)-\frac{1}{2}c_1(E).(mH) \,.$$
This is equivalent to $B^2+2B.(mH)-\Delta(E) \geq 0$.

Combined with the Hodge index theorem
$\frac{(B.H)^2}{H^2}\geq B^2$ this yields
\begin{eqnarray}\label{r34}
(B.H)^2+2mH^2\cdot (B.H)-\Delta(E)\cdot H^2& \geq &0\,.
\end{eqnarray}
Our second assumption on $m$ implies that the quadratic equation
$x^2+2mH^2\cdot x-\Delta(E)\cdot H^2 = 0$ for the
indeterminant $x$ has a positive discriminant.
Therefore, it results from (\ref{r33}) and (\ref{r34}) that  
$$B.H \geq -mH^2+\sqrt{(mH^2)^2+\Delta(E) \cdot H^2} \,.$$

By the assumption $m>\frac{a}{2H^2}-\frac{\Delta(E)}{2a}$
we eventually obtain from the last inequality $B.H>-a$.
By the very definition of $a$ this shows that $B.H$ is non negative.
\qed

\tip
{\bf Remark:}
Of course the number $a$
in theorem \ref{res5} can always be set to one.
However, the larger $a$ the sharper becomes the estimation.
In particular, if $H$ itself is the $k$th multiple of a divisor class,
then we can set $a=k$.

\section{The higher dimensional case}
Before we generalize theorem \ref{res1} to varieties of dimension at least
three,
we present a lemma which is needed in the proof of \ref{res6}.
First let us fix notations.
Let $X$ be smooth variety of dimension $n>2$ defined over an
algebraically closed field.
Furthermore, let $\Ocal_X(H)$ be a very ample line bundle.

A torsion free $X$-sheaf $F$ of rank one is a line bundle outside a set of
codimension two.
Thus, the first Chern class $c_1(L)$ is well defined.

\begin{lemma}\label{h0bound}
Let $L$ be a torsion free $X$-sheaf of rank one.
If $c_1(L).H^{n-1} < d_1$,
then there is an integer  $d_2$
only depending on $X$, $H$, and $d_1$ such that
$h^0(L)<d_2$.
\end{lemma}
\proof
We show this by induction on the dimension of $X$.
The case $n=0$ being trivial.
Suppose now that the result holds for schemes of dimension $n-1$.
Since $L(-d_1H)$ has by assumption no global sections.
Take a smooth divisor $D$ in the linear system $|d_1H|$ such that the
restriction $L|_D$ is torsion free,
and the sequence
$0 \to L(-d_1H) \to L \to L|_D \to 0$ is exact.
Then we have $h^0(L) \leq h^0(L|_D)$.
Thus, we have reduced the case to dimension $n-1$.
\qed
\begin{theorem}\label{res6}
Let $E$ be a rank two $X$-vector bundle which is semistable with respect
to $H$.
Let $l$ be an integer satisfying
$l \geq \log_2\left( \sqrt{ \left[ \frac{-\Delta(E).H^{n-2}}{H^n}
\right]_+ } +1 \right)$.
Then the restriction of $E$ to a general divisor of the linear system
$|2^lH|$ is semistable.
\end{theorem}
\proof
We follow  step by step the proof of theorem \ref{res1}.
There is no problem in generalizing most steps.
By definition of stability,
we need only the first terms of Hilbert polynomials.
Thus,
for computations with Chern classes and Hilbert polynomials
we are allowed to restrict to surfaces $S \subset X$
where $S$ is the intersection of $n-2$ divisors
of the linear system $|H|$.

In step 2 we replace $P_a$ by polynomials
$P(k)=\frac{1}{2}a_0(E) {k+n-1 \choose n-1}+a{k+n-2 \choose n-2}+ \ldots$
with 
$a< \frac{a_1(E)}{2}$.
The first Chern class of these destabilizing quotients is bounded above.
Furthermore, for $m\geq m_0 \gg 0$ we have $H^q(E(m))=0$, for all $q>0$.
Since a destabilizing quotient $Q$ on a divisor is a quotient sheaf of $E$,
it results that $H^q(Q(m))=0$, for $q>0$ and $m \geq m_0$.
It follows from lemma \ref{h0bound},
that there are upper bounds for $h^0(Q(m))$.
Thus there exists only a finite number of possible Hilbert polynomials,
for destabilizing quotients.

Now taking the minimal polynomial $P$ such that the associated Quot scheme
dominates $|mH|$ we can copy the above proof.
\qed

{\small Georg Hein,
Humboldt-Universit\"at zu Berlin,
Institut f\"ur Mathematik,
Burgstr. 26,
10099 Berlin (Germany),
hein@mathematik.hu-berlin.de
}

\end{document}